\documentclass{amsart}
\usepackage{amssymb}

\begin{document}
\def\gg{\gamma}
\font\pbglie=eufm10 
\def\xxe{{\text{\pbglie e}}}
\def\xxi{{\text{\pbglie i}}}
\newtheorem{theorem}{Theorem}[section]
\newtheorem{lemma}[theorem]{Lemma}
\newtheorem{remark}[theorem]{Remark}
\newtheorem{definition}[theorem]{Definition}
\newtheorem{corollary}[theorem]{Corollary}
\newtheorem{example}[theorem]{Example}
\font\MBOLD=cmmib10
\def\TAU{\hbox{\MBOLD\char"1C}}
\def\KAPPA{\hbox{\MBOLD\char"14}}
\def\mapright#1{\smash{\mathop{\longrightarrow}\limits\sp{#1}}}
\def\ip#1{(#1)}
\def\id{\operatorname{Id}}
\def\ffrac#1#2{{\textstyle\frac{#1}{#2}}}
\def\qedbox{\hbox{$\rlap{$\sqcap$}\sqcup$}}
\makeatletter
 \renewcommand{\theequation}{%
 \thesection.\alph{equation}}
 \@addtoreset{equation}{section}
\makeatother
\title[Riemannian V-submersion]
{Eigenforms of the Laplacian\\ for Riemannian V-submersions}
\author{Peter B. Gilkey${}^1$, Hong-Jong Kim, and JeongHyeong Park${}^2$}
\thanks{${}^1$ Research partially supported by the Max Planck Institute (Leipzig)}
\thanks{${}^2$ Research supported by R04-2000-000-00002-0 from
the Basic Research Program of the Korea Science and Engineering
Foundation}
\begin{address}{Mathematics Department, University of Oregon, Eugene Or 97403 USA}
\end{address}
\begin{email}{gilkey@darkwing.uoregon.edu}\end{email}
\begin{address}{Department of Mathematical Sciences, Seoul National University, Seoul 151-742 Korea}
\end{address}
\begin{email}{hongjong@math.snu.ac.kr}\end{email}
\begin{address}{Department of Computer and Applied Mathematics, Honam University,
 Gwangju 506-714 Korea}\end{address}
\begin{email}{jhpark@honam.ac.kr}\end{email}
\begin{abstract} Let $\pi:Z\rightarrow Y$ be a Riemannian $V$-submersion of compact $V$-manifolds.
We study when the pull-back of an eigenform of the Laplacian on
$Y$ is an eigenform of the Laplacian on $Z$, and when the
associated eigenvalue can change.
\end{abstract}
 \maketitle

\def\pext{\operatorname{ext}}
\def\pint{\operatorname{int}}

\section{introduction}\label{sect1}

We first review the situation in the smooth context. Let
$\pi:Z\rightarrow Y$ be a Riemannian submersion where $Y$ and $Z$
are compact {\bf smooth} manifolds without boundary. Let $\Phi_p$
be an eigen $p$-form of the Laplacian on $Y$ with eigenvalue
$\lambda$. Suppose that $\pi^*\Phi_p$ is an eigen $p$-form of the
Laplacian on $Z$ with corresponding eigenvalue $\mu$. This is, of
course, a rather rare phenomena that greatly restricts the
admissible geometry. We showed \cite{refGLPc,refGPa} that
$\lambda\le\mu$. If $p=0$,
then in fact the eigenvalue does not change, i.e. $\lambda=\mu$.
For $p\ge2$, we constructed examples in which the eigenvalue
actually changes, i.e. where $\lambda<\mu$. The case $p=1$
is still open.

In the present paper, we shall generalize these results to the
case where $Y$ and $Z$ are Riemannian $V$-manifolds; we must deal
with the complications introduced by the singular sets to extend
the results known in the smooth setting to the situation at hand. Throughout
the paper, we shall only deal with compact $V$-manifolds without boundary. We use the
Friedrichs extension to define the $p$-form valued Laplacian $\Delta_M^p$ on a $V$-manifold
$M$; let $E_\lambda^p(M)$ be the associated eigenspaces. We shall always suppose that the
singular set has codimension at least
$2$; this has some important analytic consequences as we shall see
in Section \ref{sect2}.

There is a rather elegant geometric
characterization of the situation when the pull back of every
eigen $p$-form on $Y$ is an eigen $p$-form on $Z$; necessarily the
eigenvalues do not change in this setting:

\begin{theorem}\label{thm1.1}
Let $\pi:Z\rightarrow Y$ be a Riemannian V-submersion of closed
$V$-manifolds where the singular sets of $Z$ and $Y$ are
of codimension at least 2.
\begin{enumerate}
\item Let $p=0$. The following conditions are
equivalent:\begin{enumerate}
\item
 $\Delta^0_Z\pi^*=\pi^*\Delta^0_Y$.
\item $\forall\lambda\ge0$, $\pi^*E^0_\lambda(Y)\subset E^0_\lambda(Z)$.
\item $\forall\lambda\ge0$, $\exists\mu(\lambda)\ge0$ such that
$\pi^*E^0_\lambda(Y)\subset E^0_{\mu(\lambda)}(Z)$.
\item The
fibers of $\pi$ are minimal.
\end{enumerate}
\item Let $p>0$. The following conditions are
equivalent:\begin{enumerate}
\item
 $\Delta^p_Z\pi^*=\pi^*\Delta^p_Y$.
\item $\forall\lambda\ge0$,
$\pi^*E^p_\lambda(Y)\subset E^p_\lambda(Z)$.
\item $\forall\lambda\ge0$, $\exists\mu(\lambda)\ge0$ such that
$\pi^*E^p_\lambda(Y)\subset E^p_{\mu(\lambda)}(Z)$.
\item The
fibers of $\pi$ are minimal and the horizontal distribution is
integrable.
\end{enumerate}\end{enumerate}
\end{theorem}

Theorem \ref{thm1.1} shows that if all the eigenspaces are
preserved, then all the eigenvalues are preserved as well. We now focus on
what happens if just a single eigenform is preserved. The
eigenvalue can not change if $p=0$; more generally, the eigenvalue
can not decrease if $p>0$. We remark that this fails in the context of manifolds with
boundary; Neumann eigenvalues can decrease \cite{Pa00}.

\begin{theorem}\label{thm1.2}
Let $\pi:Z\rightarrow Y$ be a Riemannian V-submersion of closed
$V$-manifolds where the singular sets of $Z$ and $Y$ are
of codimension at least 2.
\begin{enumerate}
\item If $0\ne\Phi\in
E_\lambda^0(Y)$ and if
 $\pi^*\Phi\in E_\mu^0(Z)$, then
$\lambda=\mu$.
\item Let $p>0$. If $0\ne\Phi\in
E_\lambda^p(Y)$ and if
 $\pi^*\Phi\in E_\mu^p(Z)$, then
$\lambda\le\mu$.
\end{enumerate}
\end{theorem}

Theorem \ref{thm1.2} is sharp if $p\ge2$. We refer to
\cite{refGLPc,refGPa} for the proof of the following result in the
smooth setting; the result in the more general context is then
immediate.

\begin{theorem}\label{thm1.3}
Let $p\ge2$ and let $0\le\lambda<\mu<\infty$ be given. There
exists a Riemannian V-submersion $\pi:Z\rightarrow Y$ and there
exists $0\ne\Phi\in E_\lambda^p(Y)$ so that $\pi^*\Phi\in
E_\mu^p(Z)$.
\end{theorem}

Here is a brief outline of the paper. In Section \ref{sect2}, we
review the definition of a $V$-manifold, the Friedrichs extension
of the Laplacian, and a basic regularity result. In Section
\ref{Sect3}, we recall some useful formulae intertwining the
coderivative on the base and on the total space of a Riemannian
submersion. We also discuss submersions in the context of
$V$-manifolds. In Section \ref{Sect4} we introduce the Hopf
fibration. We then take the quotient of the Hopf fibration by
suitably chosen cyclic group actions to construct a useful family
of $V$-submersions. We conclude the paper in Section \ref{sect-5} by
completing the proof of Theorems \ref{thm1.1} and \ref{thm1.2}.

It is a pleasant task to acknowledge useful conversations with Prof. Yorozu concerning this
paper.

\section{$V$-manifolds and $V$-submersions}\label{sect2}

The notion of a $V$-manifold was introduced by Satake
\cite{refSa}; he used the symbol `$V$' to indicate that one was
dealing with a cone-like singularity. Such manifolds are also
called {\it orbifolds}, see, for example, \cite{refSh,refT}. Their spectral geometry has
been studied by many authors; for example, see \cite{refGo}. They also appear in the study of
foliations \cite{refHe,refY}.

 In this paper, we follow the
notational conventions of
\cite{refK,refSa,refSb}. Let $O(m)$ be the orthogonal group. Let
$B_\delta$ be the ball of radius $\delta$ in $\mathbb{R}^m$
centered at the origin. If $G$ is a finite subgroup of $O(m)$,
then $G$ acts by isometries on $B_\delta$; let $B_\delta/G$
be the associated quotient space. 

It is worth noting for further
use the following fact. Let $G$ be a finite group acting on an open neighborhood $\mathcal{O}$
of the origin in $\mathbb{R}^m$. If
$G$ preserves some Riemannian metric on $\mathcal{O}$ and if $G$ fixes the origin, then in
geodesic coordinates, the action of $G$ is linerizable in the sense given above.

Let $M$ be a compact metric space. We say that $M$ is a {\it
$V$-manifold} if every point $P\in M$ has an open neighborhood
$U_P$ which is homeomorphic to $B_{\delta(P)}/G_P$ for some $\delta(P)>0$
and some finite subgroup $G_P\subset O(m)$. Let $\tilde
U_P=B_{\delta(P)}$ and let
$$\rho_P:\tilde U_P\rightarrow \tilde U_P/G_P=U_P$$
be the natural projection. Let
$$\tilde S_P:=\{\tilde Q\in\tilde U_P:\exists\gg\in G_P:\gamma\ne\id\text{ and }\gg\tilde Q=\tilde
Q\}\,.$$
be the {\it fixed point set} of $G_P$. We then have that
$$\rho_P : \tilde U_P - {\tilde S_P} \rightarrow U_P-\rho_P(\tilde
S_P)$$ is a covering projection. The {\it singular set} of $M$ is
defined to be
$$S_M:=\cup_{P\in M}\left\{\rho_P\tilde S_P\right\}\,.$$
Note that $\tilde S_P$ is the union of a finite number of linear
subspaces of $\tilde U_P$. We suppose that these subspaces all
have codimension at least $2$. We shall suppose that $M-S_M$ has the structure of a smooth
manifold and that the maps $\rho_P$ from $\tilde U_P-\tilde S_P$ to $U_P-\rho_P(\tilde
S_P)$ are local diffeomorphisms. The $V$-manifold is a
smooth manifold if $S_M$ is empty or, equivalently, if $G_P=\{I\}$
for every $P\in M$. 

We assume that the metric on $M$ is induced from a Riemannian
metric on $M-S_M$. We also assume that there is a Riemannian metric
$\tilde g_P$ on each $\tilde U_P$ which is invariant under the
action of the group $G_P$ such that $\rho_P$ is a local isometry
from $\tilde U_P-\tilde S_P$ to $U_P-S_P$.

Let $M$ be a Riemannian $V$-manifold. Let $dx$ be the associated
Riemannian measure on $M-S_M$. We shall take the Friedrichs
extension of the Laplacian from $M-S_M$ to define the Laplacian
$\Delta^p_M$ on $L^2(\Lambda^p(M))$. Let $C^\infty_0
(\Lambda^p(M-S_M))$ be the space of smooth $p$-forms which are
compactly supported in $M-S_M$. Then the $L^2$ space
$L^2(\Lambda^p(M))$ and the Sobolev space $H_1(\Lambda^p(M))$ are
defined as the completion of $C_0^\infty(\Lambda^p(M-S_M))$ with
respect to the inner products
\begin{eqnarray*}
&&(\phi,\psi)_0:=\textstyle\int_{M-S_P(M)}(\phi,\psi)dx\quad\text{and}\\
&&(\phi,\psi)_1:=\textstyle\int_{M-S_P(M)}
 \{(d\phi,d\psi)+(\delta\phi,\delta\psi)+(\phi,\psi)\}dx\,,
\end{eqnarray*}
respectively. Introduce the quadratic form
$$I^p(\phi,\psi):=\textstyle\int_M\{(d\phi,d\psi)+(\delta\phi,\delta\psi)\}dx\,.$$
The {\it Friedrichs extension} $\Delta_M^p$ is then defined
\cite{refR} by the identity:
$$(\Delta_M^p\phi,\psi)_{L^2}=I^p(\phi,\psi)\quad\text{for all}\quad\phi,\psi\in
H_1(\Lambda^p(M))\,.$$
We remark that if we removed a set  $S$ of codimension at least $2$ from a smooth manifold $M$, then this
definition of the spaces $H_1(\Lambda^p(M))$ and $\Delta_M^p$ would agree with the usual definition.

We can construct the associated spectral resolution as follows.
Let
$$\lambda_1:=\inf_{0\ne\phi\in H_1(\Lambda^p(M))}\frac{I^p(\phi,\phi)}{(\phi,\phi)_0}\,.$$
The infimum is attained by a function
$\phi_1\in\operatorname{Domain}(\Delta_M^p)$ so that
$\Delta_M^p\phi_1=\lambda_1\phi_1$. The second eigenvalue
$\lambda_2$ is then defined by setting:
$$\lambda_2:=\inf_{0\ne\phi\in
H_1(\Lambda^p(M)),\ \phi\perp\phi_1}\frac{I^p(\phi,\phi)}{(\phi,\phi)_0};$$
here `$\perp$' is with respect to the $L^2$ inner product.
Again, the infimum is attained by a function $\phi_2\in\operatorname{Domain}(\Delta_M^p)$. One
proceeds in this fashion to construct a complete orthonormal basis $\{\phi_i\}_{i=1}^\infty$
for $L^2(\Lambda^p(M))$ so that $\Delta_M^p\phi_i=\lambda_i \phi_i$.
The collection $\{\lambda_i , \phi_i\}$ is called a {\it discrete
spectral resolution} of $\Delta_M^p$ and we set
$$E_{\lambda}^p(M) =\operatorname{span}_{\{\lambda=\lambda_i\}}
\{\phi_i\}\,.$$
Again, this definition coincides with the usual definition in the smooth setting. We summarize
the discussion given above in the following Lemma and refer to
\cite{Fa01, refH, refSh} for additional related materials.

\begin{lemma}\label{lem-2.1}
Let $M$ be a closed Riemannian $V$-manifold where the singular set has codimension at least
$2$. Then:
\begin{enumerate}
\item $\Delta_{M}^p$ is self-adjoint and
non-negative.
\smallskip\item There exists a discrete spectral resolution
$\{\lambda_i , \phi_i\}$ for $\Delta_{M}^p$ where $\lambda_i\rightarrow\infty$.
\smallskip\item We have a complete orthogonal direct sum decomposition
$$L^2(\Lambda^p(M))=\oplus_\lambda E_{\lambda}^p(M)\,.$$
\end{enumerate}
\end{lemma}

The following regularity result is a central one in the subject --
we shall derive it from results of Harvey and Polking
\cite{refHaPo} but there are many other proofs see, for example,
the discussion in \cite{refC,Fa01} for $p=0$. It identifies the
eigenforms on $M$ with the smooth equivariant eigenforms on the
desingularization.

\begin{theorem}\label{thm-2.2}
Let $M$ be a closed Riemannian $V$-manifold where the singular set has codimension at least
$2$. Let
$\phi\in L^2(\Lambda^p(M))$. Then the following conditions are equivalent:

\begin{enumerate}
\item $\phi\in E_\lambda^p(M)$.
\item For any $P\in M$, there exists
$\tilde\phi_P:=C^\infty(\Lambda^p(\tilde U_P))$ so that
\begin{enumerate}
\item $\tilde\phi_P|_{(\tilde U_P - \tilde S_P)} = \rho_P^*(\phi_P|_{(U_P -
S_P)})$.
\item $\tilde\phi_P\in E_\lambda^p(\tilde U_P)$.
\item ${\gg}^*{\tilde\phi_P}=\tilde\phi_P$ for any
${\gg}\in G_P$.
\end{enumerate}\end{enumerate}
\end{theorem}

Since $S_P$ has codimension at least $2$, $C_0^\infty(\tilde U_P-\tilde S_P)$ is dense in
$H_1(\tilde U_P)$. It is now immediate from the discussion given above that Condition (2)
implies Condition (1). The converse is a smoothness result which shows that the pull back
extends smoothly across the singular set. Before establishing this
implication, we first recall a
technical Lemma.

\begin{lemma}\label{lem-2.3}
Let $\Omega$ be an open subset of $\mathbb{R}^m$ and let $A$ be a
closed subset of $\Omega$. Let $P(x,D)$ be a vector valued partial differential operator
on
$\Omega$.
\begin{enumerate}\item Assume that
$\nu:=m-2\cdot\operatorname{order}(P)\ge0$ and that the lower
Minkowski content $M_\nu(K)$ is finite for every compact set
$K\subset A$. If $\phi\in L^2_{\operatorname{loc}}(\Omega)$ and if
$P\phi=0$ on $\Omega-A$, then $P\phi=0$ on $\Omega$.
\item If $P$ is elliptic, if $\phi\in L^2_{\operatorname{loc}}(\Omega)$, and if $P\phi$ is smooth on
$\Omega$, then $\phi$ is smooth on $\Omega$.
\end{enumerate}
\end{lemma}

\begin{proof}Assertion (1) follows from Corollary 2.4 (a) of Harvey
and Polking \cite{refHaPo} who generalized earlier work of
Littman \cite{Li67}. Assertion (2) is a standard elliptic
regularity result, see, for example, Gilkey
\cite{Gi95}.\end{proof}

\begin{proof}[Proof of Theorem \ref{thm-2.2}] Let $\phi$ be an $L^2$ eigenform of the
Friedrichs extension of $\Delta$ corresponding to the eigenvalue
$\lambda$; we omit $p$ from the notation. Note that
$$||\phi||_1^2:=||\phi||_{L^2}^2+||(d+\delta)\phi||_{L^2}^2\,.$$

Let $\rho_P:\tilde U_P\rightarrow U_P=\tilde U_P/G_P$ be a local
desingularization. Let $\tilde d$ and $\tilde\delta$
denote the exterior derivative and co-derivative on $\tilde U_P$,
respectively. We then have $\rho_P^*d=\tilde d\rho_P^*$ and
$\rho_P^*\delta=\tilde\delta\rho_P^*$ since $\rho_P$ is an
isometry off the singular set. We set
$$\tilde\phi_P:=\rho_P^*\phi\,.$$
As the singularity set has codimension at least $2$ and as $\phi\in H_1(M)$,
$\tilde\phi_P\in H_1(\tilde U_P)$. The equivariance property is immediate since
$\gamma^*\rho_P^*=\rho_P^*$. To complete the proof, we must show that
$\tilde\phi_P$ extends smoothly to all of
$\tilde U_P$.

We decompose $\tilde\Delta-\lambda$ as the product of two first order operators:
$$\tilde\Delta-\lambda=(\tilde d+\tilde\delta+\sqrt\lambda)(\tilde d+\tilde\delta-\sqrt\lambda)\,.$$
We set
$$\psi:=(d+\delta-\sqrt\lambda)\phi\quad\text{and}\quad
\tilde\psi_P:=\rho_P^*\psi=(\tilde
d+\tilde\delta-\sqrt\lambda)\tilde\phi_P\,.$$ We then have that
$\tilde\psi_P$ is in $L^2$. We may express:
\begin{eqnarray*}
(\tilde d+\tilde\delta+\sqrt\lambda)\tilde\psi_P&=&(\tilde
d+\tilde\delta+\sqrt\lambda)
 (\tilde d+\tilde\delta-\sqrt\lambda)\tilde\phi_P
=(\tilde\Delta-\lambda)\tilde\phi_P\\
&=&\rho_P^*(\Delta-\lambda)\phi=0
 \qquad\text{on}\qquad\tilde U_P-\tilde S_P\,.
\end{eqnarray*}

By assumption, $\tilde {S}_P$ is the finite union of finite number
of linear subspaces of codimension at least 2 intersected with
$\tilde {U}_P$. Thus the $m-2$ dimensional upper Minkowski measure
of any compact subset of $\tilde S_P$ is finite. Lemma
\ref{lem-2.3} (1) shows:
$$(\tilde d+\tilde\delta+\sqrt\lambda)\tilde\psi_P=0\quad\text{on}\quad\tilde U_P\,.$$
Since $\tilde d+\tilde\delta+\sqrt\lambda$ is elliptic, Lemma
\ref{lem-2.3} (2) implies $\tilde\psi_P$ is smooth on $\tilde U_P$.

Since $\tilde\phi_P$ is in $H_1(\tilde U_P)$,
 $(\tilde d+\tilde
\delta-\sqrt\lambda)\tilde\phi_P=\tilde\psi_P$ in $L^2(\tilde U_P)$. Since $\tilde\psi_P$ is smooth on
$\tilde U_P$ and since this operator is elliptic, another
application of Lemma \ref{lem-2.3} (2) yields that $\tilde\phi_P$
is smooth on $\tilde U_P$ as desired.
\end{proof}
 Theorem \ref{thm-2.2} shows the pull-back eigenforms of the
Friedrichs Laplacian on $M$ are ordinary eigenforms of the
Laplacian on $\tilde U_P$ for any $P\in M$ which are invariant
under the groups $G_P$. Conversely, if we are given a collection
of eigenforms $\tilde\phi_P$ in $E_\lambda^p(\tilde U_P)$ which
are invariant under the action of the groups $G_P$ and which patch
together, then they define an eigenform of $\Delta_M^p$ on $M$.

\section{Submersions in the context of $V$-manifolds}\label{Sect3}
We begin by reviewing some of the geometry of a Riemannian submersion in the smooth setting. Let
$\pi:Z\rightarrow Y$ be a Riemannian submersion of closed smooth manifolds.
 For $z \in Z$, we decompose $T_zZ=V_z\oplus H_z$ where
$$V_z:=\ker (\pi_{*z})\quad\text{and}\quad H_z:=V_z^\perp$$
are the {\it vertical} and {\it horizontal} spaces, respectively; by assumption $\pi_*$ is an
isometry from $H_z$ to $T_{\pi z}Y$. We introduce the following notational conventions. Let
indices
$i$,
$j$, and
$k$ index local orthonormal frames $\{e_i\}$ and $\{e^i\}$ for the
vertical distributions and co-distributions ${\mathcal V}$ and
${\mathcal V}^*$ of $\pi$. Let indices $a$, $b$, and $c$ index
local orthonormal frames $\{f_a\}$ and $\{f^a\}$ for the
horizontal distributions and co-distributions ${\mathcal H}$ and
${\mathcal H}^*$ of $\pi$. If $\xi$ is a covector, then let $ \pext( \xi )$ and
$ \pint( \xi )$ denote {\it exterior multiplication} and the dual,
{\it interior multiplication}, respectively. We define tensors $\theta$ and $\omega$ and endomorphisms
$\mathcal{E}$ and $\Xi$ by setting:
\begin{equation}\label{eqn-3.a}
\begin{array}{ll}
\theta:=-\textstyle\sum_{i,a}g_Z([e_i,f_a],e_i)f^a,&
 \omega_{abi}:=\textstyle\frac12g_Z(e_i,[f_a,f_b]),\vphantom{\vrule height 11pt}\\
{\mathcal E}:=\textstyle\sum_{a,b,i}\omega_{abi}\pext_Z(e^i)
\pint_Z(f^a)\pint_Z(f^b),\vphantom{\vrule height 11pt}&
\Xi:=\pint_Z(\theta)+{\mathcal E}\,.\vphantom{\vrule height 11pt}\end{array}
\end{equation}

The tensor $\theta$ is the unnormalized mean curvature co-vector of the
fibers of $\pi$ and $\omega$ is the curvature of the horizontal
distribution. We say that the fibers are {\it minimal} if $\theta=0$.
We say that the horizontal distribution ${\mathcal H}$ is {\it integrable} if $\omega=0$.

The pull back $\pi^*$ is a linear map from $C^\infty(\Lambda^p(Y))$ to
$C^\infty(\Lambda^p(Z))$ which commutes with the exterior derivative,
i.e. $\pi^* d_Y=d_Z\pi^*$. However, $\pi^*$ does {\bf not} in
general commute with the coderivative. We refer to
\cite{refGPa} for the proof of the following result; what is
crucial for our present considerations is that the result in
question is purely local -- it does not rely on any compactness
considerations.

\begin{lemma}\label{lem-3.1}
Let $\pi:Z\rightarrow Y$ be a Riemannian submersion of Riemannian
manifolds. Then $\delta_Z\pi^*-\pi^*\delta_Y=\Xi\pi^*$ and
$\Delta_Z\pi^*-\pi^*\Delta_Y=(d_Z\Xi+\Xi d_Z)\pi^*$.
\end{lemma}

If $p=0$, the situation is simpler. If $\Phi$ is a $0$-form, i.e. a function,
then
\begin{eqnarray*}
&&\textstyle\sum_{a,b,i}\omega_{abi}\pext_Z(e^i)
\pint_Z(f^a)\pint_Z(f^b)\pi^*\Phi=0,\\
&&\textstyle\sum_{a,b,i}\omega_{abi}\pext_Z(e^i)
\pint_Z(f^a)\pint_Z(f^b)d_Z\pi^*\Phi=0,\quad\text{and}\\
&&\pint_Z(\theta)\pi^*\Phi=0\,.
\end{eqnarray*}
The following Corollary is now immediate:

\begin{corollary}\label{cor-3.2}
Let $\pi:Z\rightarrow Y$ be a Riemannian submersion of Riemannian
manifolds. Then
$\Delta_Z^0\pi^*-\pi^*\Delta_Y^0=\pint_Z(\theta)d_Z\pi^*$ on $C^\infty(Y)$.
\end{corollary}

We say that the horizontal distribution is integrable if $\omega=0$. We refer to
\cite{refGLPd} for the proof of the following result which relates $\theta$ to the
local volume element in this setting:

\begin{lemma}\label{lem-3.3}
 Let $X$ be the fiber of a Riemannian
submersion $\pi:Z\rightarrow Y$. Assume the horizontal
distribution of $\pi$ is integrable. There exist local coordinates $z=(x,y)$ on $Z$ so
 $\pi(x,y)=y$,
\begin{eqnarray*}
&&ds_Y^2=\textstyle\sum_{a,b}h_{ab}(y)dy^a\circ dy^b,\quad\text{and}\\
&&ds_Z^{2}=\textstyle\sum_{i,j}g_{ij}(x,y)dx^{i}\circ
dx^{j}+\sum_{ab}h_{ab}(y)dy^{a}\circ
 dy^{b}\,.
\end{eqnarray*}
If we set $g_{X}:=\det(g_{ij})^{1/2}$,
 then $\theta=-d_{Y}\ln(g_{X}).$
\end{lemma}

We now extend these notions to the singular setting. Let $X$ be a closed smooth manifold. We enlarge slightly the
notion of admissible charts and consider an action
\begin{equation}\label{eqn-3.b}
\gg\cdot(\tilde u,x)=(\gg\tilde u,\gg(\tilde u)x)\quad\text{on}\quad\tilde U_P\times X\,.
\end{equation}

\begin{definition}\label{defn-3.4}
 Let $Y$ and $Z$ be $V$-manifolds and let $X$ be a
smooth manifold. We say that $\pi:Z\rightarrow Y$ is {\it a
$V$-manifold fiber bundle with fiber $X$} if we can choose charts
$\tilde U_y/G_y$ over $Y$ and charts $\{\tilde U_y\times X\}/H_y$
over $Z$ so that we have a commutative diagram
\begin{equation}\label{eqn-3.c}
\begin{array}{lll}
\tilde U_y\times X&\mapright{\rho_y^Z}&\{\tilde U_y\times X\}/H_y\\
\downarrow\tilde\pi&\quad\circ&\qquad\quad\downarrow \pi\\
\tilde U_{y}&\mapright{\rho_y}&U_y:=\tilde U_y/G_y
\end{array}
\end{equation}
where $H_y$ is a subgroup of $G_y$, where $\tilde\pi$ and $\pi$ are projection on the first
factors, where the action of
$H_y$ on $\tilde U_y\times X$ is as discussed above in Equation (\ref{eqn-3.b}), and where the projections
$\rho_y^Z$ and $\rho_y$ are the associated quotient maps. We say that
$\pi$ is {\rm a Riemannian $V$-submersion} if additionally $\tilde\pi$ is a
Riemannian submersion. \end{definition}

We remark that the Riemannian metric on $\tilde U_y\times X$ is not in general a product
metric and that the decomposition in question is only local; in general, of course, $Z$ is
not $Y\times X$.

Let $\pi:Z\rightarrow Y$ be a Riemannian $V$-submersion. Let $Y_r:=Y-S_Y$ and let
$Z_r:=Z-\pi^{-1}\{S_Y\}$; these are open Riemannian manifolds. Although $Y_r$ is the regular set of $Y$, the
regular set of $Z$ may be larger than $Y_r$. Let
$\pi_r$ be the induced Riemannian submersion from $Z_r$ to $Y_r$.

\begin{lemma}\label{lem-3.5}
Let $\pi:Z\rightarrow Y$ be a Riemannian $V$-submersion of closed $V$-manifolds where
the singular sets of $Z$ and $Y$ are codimension at least 2. If $\pi_r^*\Phi\in
E_\mu^p(Z_r)$, then $\pi^*\Phi\in E_\mu^p(Z)$.
\end{lemma}

\begin{proof} Let $y\in Y$. Let $\tilde U_y/G_y$ and
$\{\tilde U_y\times X\}/H_y$ be desingularizing local charts on $Y$ and on $Z$, respectively.
By Theorem
\ref{thm-2.2},
$\tilde\Phi_y:=\rho_y^*\Phi$ extends smoothly from $\tilde U_y-\tilde S_y$ to $\tilde U_y$
with $\Delta_{\tilde U_y}^p\tilde\Phi_y=\lambda\tilde\Phi_y$. Pulling back then yields
$\tilde\pi^*\tilde\Phi_y$ is smooth on $\tilde U_y\times X$. Since
$$\Delta_{\tilde U_y\times X}^p\tilde\pi^*\tilde\Phi_y=\mu\pi^*\tilde\Phi_y\quad\text{on}\quad
\{\tilde U_y-\tilde S_y\}\times X$$
and since $\tilde S_y\times X$ has co-dimension $2$ in $\tilde U_y\times X$, we have by
continuity that:
$$\Delta_{\tilde U_y\times X}^p\tilde\pi^*\tilde\Phi_y=\mu\pi^*\tilde\Phi_y\quad\text{on}\quad
\tilde U_y\times X\,.$$
As the equivariance property is immediate, Theorem
\ref{thm-2.2} implies $\pi^*\Phi\in E_\mu^p(Z)$.
\end{proof}

We shall need to generalize the {\it fiber product} from the smooth to the singular
setting. Let
$\pi_i:Z_i\rightarrow Y$ be Riemannian $V$-submersions with fibers $X_i$ for $i=1,2$. The fiber
product is defined by setting:
$$Z=Z(Z_1,Z_2):=\{z=(z_1,z_2)\in Z_1\times Z_2:\pi_1(z_1)=\pi_2(z_2)\}\,.$$

We now study the local geometry. Choose charts on $Y$, associated charts on $Z_i$, and local projections as given
above. We shall assume that the associated
groups on $Z_i$ are the same, i.e. $H_{y,1}=H_{y,2}$ for all $y\in
Y$; this is a crucial point. The assumption $H_{y,1}=H_{y,2}$ causes no difficulty as we shall
be taking
$Z_1=Z_2$ subsequently. Then local charts for the fiber
product and the fiber product action which generalize those given in Equation (\ref{eqn-3.b}) are defined by
taking the action:
$$
\gamma\cdot(\tilde y, x_1, x_{2})=(\gamma\tilde y,
\gamma(\tilde y)x_{1}, \gamma(\tilde y)x_{2})\text{ on }
\tilde {U}_y\times X_1 \times X_2\text{ for }\gamma\in H_{y,1}=H_{y,2}\,.
$$
This shows that $\pi:Z\rightarrow Y$ is a $V$-manifold fiber bundle with fiber $X_1\times X_2$.
A similar argument can be used to show that $\pi$ is a Riemannian $V$-submersion, where a suitable
rescaling of the metric on the horizontal distribution is used, exactly as was done
in the non-singular setting \cite{refGLPc}. Let $\sigma_1(z_1,z_2):=z_1$ and
$\sigma_2(z_1,z_2):= z_2$ define maps
$\sigma_i:Z\rightarrow Z_i$. These are Riemannian $V$-submersions as well.

\begin{lemma}\label{lem-3.6} 
Adopt the notation given above. If $\Phi\in E_\lambda^p(Y)$, if
$\pi_1^*\Phi\in E_{\lambda+\varepsilon_1}^p(Z_1)$, and if $\pi_2^*\Phi\in
E_{\lambda+\varepsilon_2}^p(Z_2)$, then
$\pi^*\Phi\in E_{\lambda+\varepsilon_1+\varepsilon_2}^p(Z)$.
\end{lemma}

\begin{proof} Off the singular locus, we use the computations given in \cite{refGLPc} to see
$$
\theta_Z=\sigma_1^*\theta_{Z_1}+\sigma_2^*\theta_{Z_2}\quad\text{and}\quad
{\mathcal E}_Z\pi^*=\sigma_1^*{\mathcal
E}_{Z_1}\pi_1^*+\sigma_2^*{\mathcal E}_{Z_2}\pi_2^*\,.
$$
A straight forward application of Lemma \ref{lem-3.1} then shows
$$\Delta_Z^p\pi^*\Phi=(\lambda+\varepsilon_1+\varepsilon_2)\pi^*\Phi\quad\text{on}\quad
Z_r\,.$$ 
The desired conclusion now follows from Lemma \ref{lem-3.5}.
\end{proof}

\section{The Hopf fibration}\label{Sect4}
There is a useful family of examples which we can describe as
follows. We give the unit sphere $S^n$ in $\mathbb{R}^{n+1}$ the
standard metric $g_n$ of constant sectional curvature $+1$. We
identify $\mathbb{R}^4=\mathbb{C}^2$ to regard
$S^3\subset\mathbb{C}^2$ and we identify
$\mathbb{R}^3=\mathbb{C}\oplus\mathbb{R}$ to regard
$S^2\subset\mathbb{C}\oplus\mathbb{R}$. The {\it Hopf fibration}
$\tilde h:S^3\rightarrow S^2$ is then defined by setting
$$\tilde h (z_1,z_2):=(2z_1\bar z_2,|z_1|^2-|z_2|^2)\,.$$
One then has that $\tilde h:(S^3,g_3)\rightarrow(S^2,\ffrac14g_2)$ is a
Riemannian submersion; $(S^2,\ffrac14g_2)$ is, of course, just the
sphere of radius $\frac12$ in $\mathbb{R}^3$.

\begin{example}\label{exm4.1}
\rm Let $\mathbb{Z}_n$ be the group of $n^{\operatorname{th}}$
roots of unity in $\mathbb{C}$. We define actions $\rho_2$ and
$\rho_3$ of $\mathbb{Z}_n$ on $S^2$ and on $S^3$, respectively, by
setting
$$\rho_2(\gamma)(w,t):=(\gamma w,t)\quad\text{and}\quad\rho_3(\gamma)(z_1,z_2)=(\gamma
z_1,z_2)\quad\text{for}\quad\gamma\in\mathbb{Z}_n\,.$$ Since the actions in question are by
isometries, this gives the quotient spaces 
$$M^2_n:=S^2/\rho_2(\mathbb{Z}_n)\quad\text{and}\quad
M^3_n:=S^3/\rho_3(\mathbb{Z}_n)$$ the structure of Riemannian $V$-manifolds.
Furthermore, as the group actions are compatible with the Hopf fibration, we
have a commutative diagram
\begin{equation}\label{eqn4.a}
\begin{array}{lcl}
S^3&\mapright{\pi_3}&M^3_n\\
\downarrow \tilde h&\circ&\downarrow h\\
S^2&\mapright{\pi_2}&M^2_n
\end{array}\end{equation}
where the induced Hopf map $h$ is a
Riemannian $V$-submersion.

Let $N=(0,1)$ and $S=(0,-1)$ be the
north and south poles of $S^2\subset\mathbb{C}\oplus\mathbb{R}$.
Then $\tilde h^{-1}(N)$ is the circle $(z_1,0)$ while $\tilde h^{-1}(S)$ is the
circle $(0,z_2)$. The action of $\mathbb{Z}_n$ on $h^{-1}(N)$ is
without fixed points; thus the image of this circle in
$M^3_n$ is non-singular; the singular set of
$M^3_n$ is $\tilde h^{-1}(S)$. Thus this
example illustrates that the singular set of the total space is
not simply the inverse image of the singular set of the base. Note that the singular
sets of $M^3_n$ and of $M^2_n$ have codimension $2$.

Let $\tilde\nu_2$ be the volume form on $(S^2,\frac14g_2)$. Since $\tilde\nu_2$ is invariant
under the action of
$\mathbb{Z}_n$, it descends, by Theorem \ref{thm-2.2}, to define a harmonic $2$-form $\nu_2\in
E_0^2(M^2_n)$. One computes that $\tilde h^*\tilde\nu_2\in E_4^2(S^3)$, see
\cite{refGLP96} for details. Thus similarly
$h^*\nu_2\in E_4^2(M^3_n)$. This illustrates Theorem \ref{thm1.3} by
providing an example where the eigenvalue changes.
\end{example}

\begin{example}\label{exm4.2}
\rm Let $(p,q)$ be coprime integers and let $n=pq$. Choose
integers $a$ and $b$ so that $ap-bq=1$. We define actions $\rho_2$
and $\rho_3$ of $\mathbb{Z}_n$ on $S^3$ and on $S^2$,
respectively, by setting,
$$\rho_{2,p,q}(\gamma)(w,t):=(\gamma
w,t)\text{ and }\rho_{3,p,q}(\gamma)(z_1,z_2):=(\gamma^{ap}z_1,\gamma^{bq}z_2)
\text{ for }\gamma\in\mathbb{Z}_n\,.$$
Since the actions in question are compatible with the Hopf map, we
once again have a commutative diagram of the form given above in
Display (\ref{eqn4.a}). The action of $\mathbb{Z}_n$ on the fiber
circles $\tilde h^{-1}(N)=(z_1,0)$ and $\tilde h^{-1}(S)=(0,z_2)$ in $S^3$ is
not faithful. The isotropy group for these fiber
circles is $\mathbb{Z}_q$ and $\mathbb{Z}_p$, respectively. This
example illustrates that
the structure groups can be different over different components of
the singular set.
\end{example}

\begin{example}\label{exm4.3}
\rm Let $N$ be an arbitrary closed Riemannian manifold. We extend
the actions of $\mathbb{Z}_n$ defined in Example \ref{exm4.2} to
act trivially on $N$ to define a commutative diagram
\begin{equation}\label{eqn4.b}
\begin{array}{lcl}
S^3\times N&\mapright{\pi_3\times\id}&M^3_n\times N\\
\downarrow\tilde h\times\id&\circ&\downarrow h\times\id\\
S^2\times
N&\mapright{\pi_2\times\id}&M^2_n\times N\,.
\end{array}\end{equation}
By replacing $\nu_2$ by $\nu_2\wedge\mu_{p-2}$ for a suitably chosen eigen-form $\mu_{p-2}$ on
$N$, and by rescaling the metrics appropriately, a family of examples can be constructed
illustrating Theorem 1.3 in full generality. We omit details in the interests of brevity and
refer to \cite{refGLPc,refGPa} for further details.
\end{example}

\begin{example}\rm One can also consider the higher dimensional Hopf fibration\break $\tilde
h:S^{2n+1}\rightarrow\mathbb{CP}^n$ where the sphere $S^{2n+1}$ is given the standard
metric and where complex projective space $\mathbb{CP}^n$ is given a suitably scaled
Fubini-Study metric. Let $\mu_2\in E_0^2(\mathbb{CP}^n)$ be the Kaehler form; this restricts
to a multiple of the volume form on $S^2=\mathbb{CP}^1$. If $1\le p\le n$, then:
$$\mu_2^p\in E_0^{2p}(\mathbb{CP}^n)\quad\text{and}\quad
\tilde h^*\mu_2^p\in E_{4p(n+1-p)}^{2p}(S^{2n+1})\,.$$
We refer to Remark 3.6
\cite{refGLP96} for further details. Again, suitable cyclic group actions yield appropriate $V$-manifold
examples. 
\end{example}

\section{Proof of theorems}\label{sect-5}

\begin{proof}[Proof of Theorem \ref{thm1.2} (1)] Let $\pi:Z\rightarrow Y$ be a Riemannian $V$-submersion where
$Z$ and $Y$ are closed $V$-manifolds. We assume the singular sets of both $Y$ and $Z$ have
codimension at least $2$. Let
$0\ne\Phi\in E_\lambda^0(Y)$. Assume that
$\Psi:=\pi^*\Phi\in E_\mu^0(Z)$. Let $\rho_y:\tilde U_y\rightarrow U_y$ be an $V$-manifold
chart on $Y$. We apply Theorem \ref{thm-2.2} to see $\tilde\Phi_y:=\rho_y^*\Phi$
is smooth on $\tilde U_y$. Since $\tilde\Phi_y$ is invariant under
the action of the group $G_y$, $\Phi$ is continuous on the
quotient $U_y=\tilde U_y/G_y$. Since $\Phi$ is continuous on a compact space $Y$,
we may choose $y_0\in Y$ so
$\Phi(y_0)$ is maximal. By replacing $\Phi$ by $-\Phi$ if
necessary, we may assume without loss of generality that
$\Phi(y_0)>0$.

Let $\tilde Z_{y_0}:=\tilde
U_{y_0}\times X$ and let $\tilde\Phi_{y_0}:=\rho_{y_0}^*\Phi$. By
Theorem \ref{thm-2.2},
$$\tilde\Phi_{y_0}\in E_\lambda^0(\tilde U_{y_0})\qquad\text{and}\qquad\tilde\pi^*\tilde\Phi_{y_0}\in
E_\mu^0(\tilde Z_{y_0})\,.$$ 
We apply Corollary
\ref{cor-3.2} to the Riemannian submersion
$\tilde\pi:\tilde Z_{y_0}\rightarrow\tilde U_{y_0}$ to see:
\begin{equation}\label{eqn-5.a}
(\mu-\lambda)\tilde\pi^*\tilde\Phi_{y_0}=\{\Delta_{\tilde Z}^0\pi^*-\pi^*\Delta_{\tilde
Y}^0\}\tilde\Phi_{y_0}=\pint(\tilde\theta)d_{\tilde Z}\tilde\pi^*\tilde\Phi_{y_0}
=\pint(\tilde\theta)\tilde\pi^*d_{\tilde Y}\tilde\Phi_{y_0}\,.
\end{equation}
Choose $\tilde z_0$ so $\tilde\pi\tilde z_0=\tilde y_0$. Since $\tilde\Phi_{y_0}$ has a maximum at
$\tilde y_0$,
$(\tilde\pi^*d_{\tilde Y}\tilde\Phi_{y_0})(\tilde z_0)=0$. Since
$(\pi^*\tilde\Phi_{y_0})(\tilde z_0)>0$, evaluating Equation
(\ref{eqn-5.a}) at $\tilde z_0$ implies $\mu=\lambda$.
\end{proof}

\begin{proof}[Proof of Theorem \ref{thm1.2} (2)] Let $\pi:Z\rightarrow Y$ be a Riemannian $V$-submersion where
$Z$ and $Y$ are closed $V$-manifolds. We assume the singular sets of both $Y$ and $Z$ have
codimension at least
$2$. Let
$p>0$, and let
$0\ne\Phi\in E_\lambda^p(Y)$. Assume that
$\pi^*\Phi\in E_\mu^p(Z)$. We wish to show $\lambda\le\mu$.

We generalize the argument given in \cite{refGLPc}. Let $Z_0:=Z$ and let $Z_1:=Z(Z_0,Z_0)$ be
the fiber product. Let $\varepsilon_0:=\mu-\lambda$. Then by Lemma \ref{lem-3.6},
$$\pi_1^*\Phi\in E_{\lambda+2\varepsilon_0}^p(Z_1)\,.$$
We now inductively set
$Z_n:=Z(Z_{n-1},Z_{n-1})$ and apply the same argument to see 
$$\pi_n^*\Phi\in
E_{\lambda+2^n\varepsilon_0}^p(Z_n)\,.$$
Since $\Delta_{Z_n}^p$ is a non-negative operator by
Theorem \ref{lem-2.1}, we have $\lambda+2^n\varepsilon_0\ge0$ for all $n$. This implies
$\varepsilon_0\ge0$ and hence $\mu\ge\lambda$.
\end{proof}

\begin{proof}[Proof of Theorem \ref{thm1.1}] We extend the arguments given in \cite{refGPa}. Let
$\pi:Z\rightarrow Y$ be a Riemannian
$V$-submersion where
$Z$ and $Y$ are closed $V$-manifolds. We assume the singular sets of both $Y$ and $Z$ have
codimension at least
$2$.

We first show that Assertion (1d) implies Assertion (1a) and that Assertion (2d) implies
Assertion (2a). Assume off the singular set that the fibers of
$\pi$ are minimal and, if $p>0$, that the horizontal distribution is integrable. Let
$0\ne\Phi\in E_\lambda^p(Y)$. Then Lemma \ref{lem-3.1} and Corollary \ref{cor-3.2} imply
that
$$\Delta_Z^p\pi^*\Phi=\lambda\pi^*\Phi\quad\text{on}\quad Z_r\,.$$
Thus by Lemma \ref{lem-3.5},
$\pi^*\Phi\in E_\lambda^p(Z)$. This shows that
$$\pi^*E_\lambda^p(Y)\subset E_\lambda^p(Z)$$
for all $\lambda$; the intertwining relations of Assertions (1a) and (2a) now follow as the
span of the eigenspaces is dense in the appropriate topology.

It is immediate that Assertion (1a) implies Assertion (1b) and that Assertion (2a) implies Assertion (2b).
Similarly Assertion (1b) implies Assertion (1c) and Assertion (2b) implies Assertion (2c). We complete the
proof by showing that Assertion (1c) implies Assertion (1d) and that Assertion (2c) implies Assertion (2d).

Suppose that $\pi^*E_\lambda^p(Y)\subset E_{\mu(\lambda)}^p(Z)$ for all
$\lambda$. Let
$\Phi_\lambda\in E_\lambda^p(Y)$. By Lemma
\ref{lem-3.1}
\begin{equation}\label{eqn-5.b}
(\mu-\lambda)\pi^*\Phi_\lambda=\{d_Z\
(\pint_Z(\theta)+{\mathcal{E}})
 +(\pint_Z(\theta)+{\mathcal{E}})d_Z\}\pi^*\Phi_\lambda
\quad\text{on}\quad Z_r\,.
\end{equation}

Suppose first that $p=0$. By Theorem \ref{thm1.2} (1), we have $\mu(\lambda)=\lambda$. We use Corollary
\ref{cor-3.2} to rewrite Equation (\ref{eqn-5.b}) in the form:
$$
\pint_Z(\theta)\pi^*d_Y\Phi_\lambda=0\quad\text{on}\quad Z_r\,.
$$
Let $\Phi$ be any smooth function which is compactly supported near a regular point $y$ of
$Y$. We can approximate $\Phi$ in the Sobolev-$H_1$ topology as a finite
sum of eigenfunctions. Thus
$$\pint_Z(\theta)\pi^*d_Y\Phi=0\,.$$
Since $\theta$ is a horizontal
co-vector and since $\Phi(y)$ is arbitrary, this implies $\theta=0$ at any point $z\in\pi^{-1}(y)$. Thus $\theta$
vanishes on $Z_r$; passing to a local desingularization, we conclude
$\tilde\theta$ vanishes everywhere by continuity. This completes the proof of Theorem
\ref{thm1.1} (1).

We now suppose that $p>0$. Let $\pi_{\mathcal {H}}$ be
orthogonal projection from $\Lambda^p(Z_r)$ to
$\Lambda^{p}(\mathcal{H})$. Let $\Phi_\lambda\in E_{\lambda}^{p}(Y)$. We apply $(1-\pi_{\mathcal{H}})$ to
Equation (\ref{eqn-5.b}) to see:
\begin{eqnarray*}
0&=&(\mu-\lambda)\pi_{\mathcal{H}}\pi^*\Phi_\lambda\\&=&
(1-\pi_{\mathcal {H}})\{d_Z(\pint_Z(\theta)+\mathcal {E})
 +(\pint_Z(\theta)+\mathcal{E})d_Z \}\pi^*\Phi_\lambda\quad\text{on}\quad Z_r\,.
\end{eqnarray*}
The natural domain of this identity is the Sobolev space $H_1(Y_r)$ and, as the eigenfunctions are dense in
$H_1(Y_r)$, by continuity we then have
\begin{equation}\label{eqn-5.c}
0=(1-\pi_{\mathcal {H}})\{d_Z(\pint_Z(\theta)+\mathcal {E})
 +(\pint_Z(\theta)+\mathcal{E})d_Z\}\pi^*\Phi\quad\text{for}\quad\Phi\in H_1(Y_r)\,.
\end{equation}
Let $\pi z_0=y_0\in Y_r$. Choose
$F\in C_0^\infty(Y_r)$ so that $F(y_0)=0$. Let
$\xi:=dF(y_0)$. Since
$\pint_Z(\theta)+{\mathcal E}$ is a $0{}^{\text{th}}$ order
operator, we apply Equation (\ref{eqn-5.c}) to the product $F\Phi$ and evaluate
at $z_0$ to see
$$0=(1-\pi_{\mathcal H})
 \{\pext_Z(\pi^*\xi)(\pint_Z(\theta)
 +{\mathcal E})+(\pint_Z(\theta)+{\mathcal E})\pext_Z(\pi^*\xi)\}
 \pi^*{\Phi(y_0 )}\,.$$
Since
$$0=(1-\pi_{\mathcal H})
 \{\pext_Z(\pi^*\xi)\pint_Z(\theta)
 +\pint_Z(\theta)\pext_Z(\pi^*\xi)\}\pi^*,$$
and since ${\mathcal E}$ always introduces a vertical covector, we
conclude
\begin{equation}\label{eqn-5.d}
\begin{array}{l}
0=(1-\pi_{\mathcal H})\{\pext_Z(\pi^*\xi){\mathcal E}+{\mathcal
E}\pext_Z(\pi^*\xi)\}\pi^*\\
\phantom{0}=\{\pext_Z(\pi^*\xi){\mathcal E}+{\mathcal
E}\pext_Z(\pi^*\xi)\}\pi^*\,.\vphantom{\vrule height 11pt}\end{array}
\end{equation}
Let $\{f^a, e^i\}$ be the
orthonormal frames of the dual distributions ${\mathcal H}^*$ and ${\mathcal V}^*$. Adopt the notation of
Equation (\ref{eqn-3.a}). To simplify the notation, set
$$\xxe_i:=\pext_Z(e^i),\quad
\xxe_a:=\pext_Z( f^a),\quad\text{and}\quad\xxi_a:=\pint_Z( f^a)\,.$$
We then have the Clifford commutation relations:
$$\xxe_a\xxi_b+\xxi_b\xxe_a=\delta_{ab}\,.$$
Choose $F$ so that $\pi^*\xi(z_0)= f^c(z_0)$ and
apply Equation (\ref{eqn-5.d}) to compute at $z_0$ that:
\begin{eqnarray*}
0&=&\textstyle\sum_{a,b,i}\omega_{abi}\{\xxe_c\xxe_i\xxi_a\xxi_b+\xxe_i\xxi_a\xxi_b\xxe_c\}
 =\sum_{a,b,i}\omega_{abi}\xxe_i\{-\xxe_c\xxi_a\xxi_b+\xxi_a\xxi_b\xxe_c\}\\
 &=&\textstyle\sum_{a,b,i}\omega_{abi}\xxe_i\{\xxi_a\xxe_c\xxi_b+\xxi_a\xxi_b\xxe_c-\delta_{ac}\xxi_b\}\\
&=&\textstyle\sum_{a,b,i}\omega_{abi}\xxe_i\{-\xxi_a\xxi_b\xxe_c+\xxi_a\xxi_b\xxe_c-\delta_{ac}\xxi_b
 +\delta_{bc}\xxi_a\}
 =-2\textstyle\sum_{b,i}\omega_{cbi}\xxe_i\xxi_b.
\end{eqnarray*}
As $p\ge1$, we may conclude that $\omega=0$ on $Z_r$ and consequently ${\mathcal H}$ is integrable.

We must now show the fibers are minimal. Let $d_X$ denote exterior
differentiation along the fiber. We set ${\mathcal E}=0$ and use
Equation (\ref {eqn-5.c}) to compute
\begin{eqnarray*}
0&=&(1-\pi_{\mathcal {H}})\{d_Z\pint_Z(\theta)
 +\pint_Z(\theta)d_Z\}\pi^*\\
&=&(1-\pi_{\mathcal{H}})d_X\pint_Z(\theta)\pi^*\quad\text{on}\quad C_0^\infty(\Lambda^p(Y_r))\,.
\end{eqnarray*}
This implies $\theta$ is constant on the
fibers so $\theta=\pi^*\Theta$ is the pull back of a globally
defined 1-form away from the singular set on the base. 

We apply Lemma \ref{lem-3.3}. Let
$d\nu^e_x$ be the Euclidean measure and let 
$$\psi(y):=\textstyle\int_X g_X(x,y)d\nu^e_x$$ be the volume of the fiber $\pi^{-1}(y)$ for $y\in Y_r$.
Then
\begin{eqnarray*}
&&d_Y\psi(y)=d_Y\textstyle\int_X g_X(x,y)d\nu^e_x
 =\textstyle\int_X (g_X g_X^{-1}d_Yg_X)(x,y)d\nu^e_x\\
&&\quad=-\textstyle\int_X g_X(x,y)\theta(x,y)d\nu^e_x
 =-\Theta(y)\textstyle\int_X g_X(x,y)d\nu^e_x\\
&&\quad=-\Theta(y)\psi(y)\quad\text{so}\\
&&\theta=-\pi^*d_Y\ln\psi\quad\text{on}\quad Y_r\,.
\end{eqnarray*}

If $y$ is a singular point, we let $\tilde U_y$ be the desingularization of $Y$ and $\tilde
U_y\times X$ be the desingularization of $Z$. Since $\tilde{\mathcal{E}}=0$ on $(\tilde
U_y-\tilde S_y)\times X$, we have $\tilde{\mathcal{E}}=0$ on $\tilde U_y\times X$ by
continuity since the singular set has codimension at least $2$. Thus we can apply exactly the
same argument given above to see $\tilde\psi:=\rho_y^*\psi$ extends to a smooth function
$\tilde\psi_y$ on all of
$\tilde U_y$.

We define
a conformal variation of the metric on the vertical distribution
which leaves the metric on the horizontal distribution unchanged by setting
$$g(t)_Z=\psi^{2t} ds_{\mathcal V}^2+ds_{\mathcal H}^2\quad\text{on}\quad Z_r\,.$$
The argument given above shows this variation extends to
the desingularization to define a smooth $1$ parameter of Riemannian
$V$-submersions.

We have that $\pi:Z(t)\rightarrow Y$ is a Riemannian submersion with integrable horizontal
distribution. We use Lemma
\ref{lem-3.3} to see $\theta(t)=(1+t\dim(X))\theta$ away from the singular set and thus
\begin{eqnarray*}
 &&\Delta_{Z}^p\pi^*-\pi^*\Delta_Y^p
 =(1+t\dim(X))(d_Z\pint_Z(\theta)+\pint_Z(\theta)d_Z)\pi^*\\
 &=&(1+t\dim(X))(\Delta_Z^p\pi^*-\pi^*\Delta_Y^p)\,.
\end{eqnarray*}
Let $\Phi\in E_\lambda^p(Y)$ and let $\pi^*\Phi\in E_\mu^p(Z)$. Set
$\varepsilon=\mu-\lambda$. Then
$$\Delta_Z^p\pi^*\Phi=\{\lambda+(1+t\dim(X))\varepsilon\}\pi^*\Phi\quad\text{on}\quad
Z_r\,.$$
Consequently, by Lemma \ref{lem-3.5},
$$\pi^* E_{\lambda}^p(Y)
 \subset E_{\lambda+(1+t\dim(X))\epsilon(\lambda)}^p(Z(t))\,.$$
By Lemma \ref{lem-2.1}, $\lambda+(1+t\dim(X))\epsilon(\lambda)\ge0$.
Since $t$ is arbitrary, $\epsilon(\lambda)=0$. Thus
\begin{equation}\label{eqn-5.e}
(d_Z\pint_Z(\theta)+\pint_Z(\theta)d_Z)\pi^*=0\quad\text{on}\quad E_\lambda^p(Y)\,.
\end{equation}
Since these eigenspaces are dense in the Sobolev space
$H_1$, Equation (\ref {eqn-5.e}) continues to be valid on
$H_1$. Let $y_0\in Y_r$. Let $f\in\mathcal{H}^*(z_0)$ where $\pi(z_0)=y_0$. Choose
$\Phi\in C_0^{\infty}(U_y)$ so that 
$$\Phi(y_0)=0\quad\text{and}\quad\pi^*d\Phi(y_0)=f\,.$$
Let $\Psi\in C_0^\infty(Y_r)$ with $\Psi(y_0)$
arbitrary. We apply Equation (\ref {eqn-5.e}) to the product
$\Phi\Psi$ and evaluate at $y_0$ to see:
$$(\pext_Z(f)\pint_Z(\theta)+\pint_Z(\theta)\pext_Z(f))\pi^*\{\Psi(y_0)\}=0\,.$$
Since
$$\pext_Z(f)\pint_Z(\theta)+\pint_Z(\theta)\pext_Z(f)=g_Z(f,\theta),$$
$g_Z(f,\theta)(z_0)=0$. Since $\theta$ is horizontal and since $f$ was an arbitrary
horizontal covector, we conclude $\theta$ vanishes away from singular set; passing to a
local desingularization, we complete the proof by checking that $\tilde\theta=0$ everywhere by
continuity.
\end{proof}


\begin{thebibliography}{AAA}

\bibitem{refC} Y.J. Chiang, {\it Spectral geometry of V-Manifolds and its application
to Harmonic maps}, Proc. Sympos. Pure Math.
{\bf 54} (1993), 93--99.

\bibitem{Fa01} C. Farsi, {\it Orbifold spectral theory}, Rocky
Mountain J. Math. {\bf 31} (2001), 215--235.


\bibitem{Gi95} P. Gilkey, {\it Invariance theory, the heat
equation, and the Atiyah-Singer index theorem $2^{nd}$ ed.}, Studies in Advanced Mathematics,
CRC Press, Boca Raton, FL, (1995). 

\bibitem{refGLP96} P. Gilkey, J. V. Leahy, and J. H. Park, {\it The spectral geometry of
the Hopf fibration},
J. Phys. A {\bf 29} (1996), 5645--5656.

\bibitem{refGLPd} ---, {\it Spinors, spectral geometry,
and Riemannian submersions}, Lecture Notes Series {\bf 40},
Seoul National University, Research Institute of Mathematics, Global Analysis Research Center, Seoul, (1998). 

\bibitem{refGLPc} ---, {\it Eigenvalues of the form valued Laplacian
for Riemannian submersions}, Proc. Amer. Math. Soc. {\bf 126} (1998),
1845--1850.

\bibitem{refGPa} P. Gilkey and J. H. Park, {\it Riemannian submersions which preserve the eigenforms
 of the Laplacian}, Illinois J. Math, {\bf 40} (1996), 194--201.

\bibitem{refGo} C. Gordon, D. Webb, and S. Wolpert, {\it Isopectral plane domains and
surfaces via Riemannian orbifolds}, Invent. Math. {\bf 110} (1992), 1--22.

\bibitem{refHaPo} R. Harvey and J. Polking, {\it Removable Singularities of solutions of linear
partial differential equations}, Acta. Math. {\bf 125} (1970),
39--56.

\bibitem{refHe} J. Hebda, {\it An example relevant to curvature pinching theorems for
Riemannian foliations}, Proc. Amer. Math. Soc. {\bf
114} (1992), 195--199.

\bibitem{refH} C. Hodgson and J. Tysk, {\it Eigenvalue estimate and isoperimetric inequalities
for cone manifolds}, Bull. Austral. Math. Soc. {\bf 47} (1993),
127--143.

\bibitem{refK} H. Kitahara, {\it On a parametrix form in a certain V-submersion},
Springer Lecture Notes in Math. {\bf 792} (1980), 264--298.

\bibitem{Li67} W. Littman, {\it Polar sets and removable
singularities of partial differential equations}, Ark. Mat. {\bf 7} (1967) 1--9.

\bibitem{Pa00} J. H. Park, {\it The spectral geometry of Riemannian submersions
for manifolds with boundary}, Rocky Mountain J. Math. {\bf 30} (2000), 353--369.

\bibitem{refR} M. Reed and B. Simon, {\it Methods of modern mathematical physics}, vol I-IV,
Academic Press, Inc. New York, (1980).

\bibitem{refSa} I. Satake, {\it On a generalization of the notion of manifold}, Proc. Nat. Acad. Sci.
U.S.A {\bf 42} (1956), 359--363.

\bibitem{refSb} I. Satake, {\it The Gauss-Bonnet theorem for $V$-manifolds}, J. Math. Soc. Japan, {\bf 9}
(1957), 464--492.

\bibitem{refSh} T. Shioya, {\it Eigenvalues and suspension structure of compact
Riemannian orbifolds with positive Ricci curvature}, Manuscripta
Math. {\bf 99} (1999), 509--516.

\bibitem{refT} W. P. Thurston, {\it Three-dimensional geometry and topology}, Princeton Mathematical Series,
{\bf 35},
Princeton University Press, Princeton, NJ (1997).

\bibitem{refY} S. Yorozu, {\it Leaf space of a certain Hopf $r$-foliation}, Nihonkai Math. J.
{\bf 10} (1999), 117--130.

\end{thebibliography}
\end{document}